\documentclass[12pt]{article}
\usepackage{amsmath,amsfonts,amssymb,amsthm}

\def\noi{\noindent}

\def\pf{\noi{\bf Proof.\ \,}}
\def\eop{{$\square$}}

\def\em{\it }        
\def\ip#1#2{\langle #1\mid #2 \rangle}

\def\labtt#1{\label {#1}}

\def\refpp#1{(\ref {#1})}

\def\CC{{\mathbb C}}

\def\NN{{\mathbb N}}
\def\QQ{{\mathbb Q}}
\def\RR{{\mathbb R}}

\def\ZZ{{\mathbb Z}}

\def\Z{{\mathbb Z}}

\def\<{\langle}
\def\>{\rangle}

\def\vac{\hbox{\bf 1}} 

\begin{document}

\newcommand{\aut}{\mathrm{Aut}}
\newcommand{\M}{\mathbb{M}}
\newcommand{\dih}[2]{DIH_{#1}(#2)}
\newcommand{\drt}[1]{\frac{1}{\sqrt{2}} #1}
\newcommand{\drtp}[1]{  \frac{1}{\sqrt{2}}\left( #1\right)}
\newcommand{\drtpp}[1]{ \left( \frac{1}{\sqrt{2}}\left( #1\right)\right)}
\newcommand{\mspan}{\mathrm{span}}
\newcommand{\Rspan}[1]{\mathrm{span}_{\mathbb{R}}(#1)}
\newcommand{\Dih}[1]{DIH_{#1}}

\newtheorem{thm}{Theorem}[section]
\newtheorem{prop}[thm]{Proposition}
\newtheorem{lem}[thm]{Lemma}
\newtheorem{coro}[thm]{Corollary}
\newtheorem{conj}[thm]{Conjecture}

\newtheorem{hyp}[thm]{Hypothesis}

\newtheorem{rem}[thm]{\bf Remark}
\newtheorem{de}[thm]{\bf Definition}
\newtheorem{nota}[thm]{\bf Notation}
\newtheorem{ex}[thm]{\bf Example}
\newtheorem{proc}[thm]{\bf Procedure}

\begin{center} 
{\Large  \bf  Determinants for integral forms in lattice type vertex operator algebras}

\bigskip

3 April, 2018

\vskip 1cm

{\Large    }

Chongying Dong

Department of Mathematics,

University of California,

Santa Cruz, CA 95064 USA

\&

School of Mathematics and Statistics, 

Qingdao University, 

Qingdao 266071 CHINA

{\tt dong@ucsc.edu}

\vskip 0.5cm
and
\vskip 0.5cm

Robert L. Griess Jr.

Department of Mathematics,

University of
Michigan,

Ann Arbor, MI 48109-1043  USA

{\tt rlg@umich.edu}

\smallskip

\end{center}

\begin{abstract}
We prove a determinant formula for the standard integral form of a lattice vertex operator algebra.
\end{abstract}

\section{Introduction}

We have studied group-invariant integral forms in vertex operator algebras \cite{DG,DG2}.   In this article, we study standard integral forms in lattice vertex operator algebras and give the determinant of each homogeneous piece as a particular integral power of determinant of the input positive definite  even integral lattice.  When the lattice is unimodular, all these homogeneous pieces have determinant 1, already  proved in \cite{DG}.   For lattices of other determinants, there did not seem to be an obvious answer.
 Borcherds stated without proof  in \cite{B1986} that the determinant of a homogeneous piece was some (unspecified) integral power of the input lattice.

The standard integral form $V_{L.\ZZ}$  for a lattice vertex operator algebra $V_L$ is reviewed in Section 4, Definition \refpp{standardintform}.   Lemma \refpp{l1} shows that our main theorem \refpp{maintheorem} is reduced to a study of determinants for integral forms within a certain symmetric algebra.  The latter determinant is therefore our main object of study in this article.

\section{Background}

\begin{lem}\labtt{ballsinurns}   If $x_1, \dots, x_k$ are variables, then the number of monomials $x_1^{a_1}\cdots x_k^{a_k}$, $a_i \in \ZZ_{\ge 0}$, of total degree $n$ is ${n+k-1}\choose {k-1}$.
\end{lem}

\pf   This is essentially a counting result, called Balls in Urns.   Monomials correspond to the set of $k-1$ marker balls to be chosen among a set of $n+k-1$ balls arranged in a straight line.   One adds marker ball 0 at the very beginning and marker ball $k$ at the very end.    The sequence $a_1, \dots, a_k$  gives the lengths of the gaps between successive marker balls.  \eop

\begin{lem}\labtt{detindex}
If $J \le K$ are finite rank lattices and the index $|J:K|$ is finite, then $det(K)=|J:K|^2 det(J)$.
\end{lem}

\section{Symmetric algebas}

\begin{nota}\labtt{01} Let $H$ be a $k$-dimensional vector space over  $\CC$ and let $t$ be a variable.  For $r\ge 1$, let $H_r=H\otimes  \CC t^{-r}$  be a copy of $H$, defined to have degree $r$, and set $h(-r)=h\otimes t^{-r}$ for
$h\in H.$ We shall work in the symmetric algebra $M(1)=\mathbb S [H[t^{-1}]t^{-1}]=\CC[h(-r)|h\in H]$ where $$H[t^{-1}]t^{-1}=\oplus_{r\geq 1}H\otimes \CC t^{-r}.$$
  Then $$M(1)=\oplus_{n\geq 0}M(1)_n$$
  is graded such that $M(1)_n$ is spanned by $h_1(-r) \cdots h_p(-r_p)$ for $h_1,...,h_p\in H$ and $r_1,...,r_p\in \NN$ with $\sum_ir_i=n.$
For a sequence $a=(a_1, \dots, a_n, \dots )$ of nonnegative integers which is almost all zero, define $wt(a):=\sum_{j\ge 1} \ ja_j$.
For $n\ge 0$, define $\mathcal A (n)$ to be the set of such $a$ of weight $n$.  Note that $\mathcal A (0)=\emptyset$.   

Suppose that $x_1, \dots, x_k$ is a basis of $H$.   Then $L=\Z x_1+\cdots +\Z x_k$ is a free abelian group of rank $r$ in $H$ and $M(1)=\CC[x_1(-r), \cdots x_k(-r)| r \geq 1].$   Define $B(a):=B(x_1, \dots x_k; a)$ to be the $\ZZ$-span of all words $w_1\cdots w_n$ in $M(1)$ where $w_i$ is a product of length $a_i$ in the variables $x_j(-i)$, for $j\in \{1,2,\dots ,k\}$.
Finally, for an integer $n\ge 0$, define $B_L(n):=B(x_1, \dots x_k; n) :=\oplus _{a\in \mathcal A (n)} B(a)$ and $B_L:=\oplus _{n\ge 0} B_L(n)=\ZZ[x_1(-r), \cdots x_k(-r)|r \geq 1]$.  Then $B_L$ is a subring of $R.$     Note that these objects are unchanged if  $x_1, \dots, x_k$ is replaced by any basis of the $\ZZ$ span of $x_1, \dots, x_k$.
\end{nota}

\begin{nota}\labtt{01'}
Let $A_L$ be the $\ZZ$-submodule of $M(1)$ generated by $s_{\alpha,n}$ for $\alpha_i\in \{x_1,\dots x_k\}$ and $n\geq 0$ where
$$E^-(-\alpha,z):=\exp\left(\sum_{n>0}\frac{\alpha(-n)}{n}z^n\right)=\sum_{n\geq 0}s_{\alpha,n}z^n.$$

Although we do not use the vertex operator algebra structure on $M(1)$,  we  use the notations $E^-(-\alpha,z)$ and $s_{\alpha,n}$ from \cite{flm} and \cite{DG} here.
Then $A_L$ has a $\ZZ$-base $B_1B_2\cdots B_k$ where
$$B_i=\{s_{x_i,n_1}\cdots s_{x_i,n_q}|n_1\geq \cdots \geq n_q\geq 0\}$$
and $B_1\cdots B_k=\{u_1u_2\cdots u_k|u_i\in B_i\}.$
\end{nota}
We also set $A_L(n)=A_L\cap M(1)_n$ for all $n.$ The following result will be useful in computing the determinants for the lattice vertex operator algebras.

\begin{lem}\label{l1}  $B_{L}$ is a subring of $A_L$ and the index $[A_L(n):B_L(n)]$ is independent of the base $\{x_1,...,x_k\}$ for any $n\geq 0.$

\end{lem}
\pf   We first prove that $B_L$ is a subring of $A_L.$ It is good enough to show that $\alpha(-n)\in B_L$ for $\alpha\in\{x_1,...,x_k\}$ and $n\geq 0.$  Note that
$$E^-(-n\alpha,z)=E^-(-\alpha,z)^n=(\sum_{m\geq 0}s_{\alpha,m}z^m)^n.$$
So the coefficient $c_n$ of $z^n$ in $E^-(-n\alpha,z)$ lies in $A_L.$ Clearly, $c_n$ is also the coefficient of $z^n$ in $(\sum_{m=0}^ns_{\alpha,m}z^m)^n.$
A straightforward computation shows that
$a_n=\alpha(-n)+u$ where $u$ is a $\ZZ$-linear combination of elements of the form $s_{\alpha,m_1}s_{\alpha,m_2}\cdots $ with$m_i<n$ and $m_1+m_2+\cdots =n.$ As a result, $\alpha(-n)\in A_L.$

To show that the index $[A_L(n):B_L(n)]$ is independent of the base $\{x_1,...,x_k\}$ for any $n\geq 0$, we let $\{y_1,...,y_k\}$ be another basis of $H$
and $K=\ZZ y_1+\cdots +\ZZ y_k.$
Then a group isomorphism  $f$ from $L$ to $K$ by sending $x_i$ to $y_i$  induces a ring isomorphism $\hat f$ from $R$ to itself  such
that $\hat f(A_L)=A_K$ and $\hat f(B_L)=B_K.$ It is evident that $\hat f$ is a degree preserving map. As a result, $\hat f(A_L(n))=A_K(n)$ and $\hat f(B_L)=B_K.$ Thus, $[A_L(n):B_{L}(n)]=[A_K(n):B_K(n)].$
 \eop

Note that both $A_L$ and $B_L$ are  $\ZZ$-forms of $M(1).$

\begin{lem} \labtt{02}
(i) $rank(S^m(H))={{m+k-1}\choose {k-1}}$;

(ii) $rank(B(a))=\prod_{j=1}^n {{a_j+k-1}\choose {k-1}}$;

(iii) $wt(B(a))=\sum_{j\ge 1} \ ja_j$;

(iv) $rank(B(n))=\sum_{a: wt(a)=n} rank(B(a))$.
\end{lem}
\pf   Straightforward, with Lemma \refpp{ballsinurns}.
\eop
\bigskip

So far,  there is no bilinear form in this discussion.  We shall introduce forms later, after Corollary \ref{06}.

\bigskip

We now compare what happens to the $B(a)$ when $x_1, \dots, x_n$ is replaced by another basis.  We already noted in Notation \ref{01}  that $B(x_1, \dots x_k; a) = B(y_1, \dots y_k; a)$ if $span_{\ZZ}(x_1, \dots x_k)=span_{\ZZ}(y_1, \dots y_k)$.

Using the proof of Lemma \ref{l1} we can easily have:
\begin{lem} \labtt {03}  If $x_1, \dots x_k$ and $y_1, \dots y_k$ are bases and if
$span_{\ZZ}(x_1, \dots x_k)$ contains $span_{\ZZ}(y_1, \dots y_k)$, then for any invertible linear transformation $T$ on $H$,
$B(x_1, \dots x_k; a)/B(y_1, \dots y_k; a) \cong B(Tx_1, \dots Tx_k; a)/B(Ty_1, \dots Ty_k; a)$.   In particular, we have equality of indices $|B(x_1, \dots x_k; a) : B(y_1, \dots y_k; a) |=|B(Tx_1, \dots Tx_k; a) : B(Ty_1, \dots Ty_k; a)|$.
\end{lem}

\begin{lem} \labtt {04} Suppose that $p>0$ is an integer.   Then
$B(x_1,  x_2, \dots , x_k; a)$ contains $B(px_1,  px_2, \dots , px_k; a)$ and the index is
$p^{N(k,a)}$ where
$N(k,a):= \prod_{j=1}^n a_j{{a_j+k-1}\choose {k-1}}$.
\end{lem}
\pf
The free abelian group $B(x_1,  x_2, \dots , x_k; a)$ has basis consisting of monomials in the $x_t(-j)$.
Such a monomial has a unique expression $w_1\cdots w_{wt(a)}$,
where $w_j$ is a monomial in the $x_t(-j)$.
 There are ${{a_j+k-1}\choose {k-1}}$ such $w_j$.
 The formula for $N(k,a)$ is now clear.
\eop

\begin{lem} \labtt {05} Suppose that $p>0$ is an integer.   Then
$B(x_1,  x_2, \dots , x_k; a)$ contains $B(px_1,  x_2, \dots , x_k; a)$ and the index is
$p^{{\frac 1k} N(k,a)}$ where
$N(k,a):=  \prod_{j=1}^n a_j{{a_j+k-1}\choose {k-1}}$.
\end{lem}
\pf Observe that we have a chain  $$span_{\ZZ}(x_1, x_2, x_3, \dots x_k) > span_{\ZZ}(px_1, x_2, x_3, \dots x_k) > $$ $$span_{\ZZ}(px_1, px_2, x_3, \dots x_k) > \dots > span_{\ZZ}(px_1,p x_2, px_3, \dots px_k).$$  By Lemma \refpp{03},
the indices for each containment
$$B(x_1, x_2, x_3, \dots x_k; a) > B(px_1, x_2, x_3, \dots x_k; a) > $$
$$B(px_1, px_2, x_3, \dots x_k; a) > \dots > B(px_1,p x_2, px_3, \dots px_k; a)$$ are equal.   We then deduce the result from Lemma \refpp{04}.
\eop

\begin{coro}\labtt{06}
In the notation of Lemma \refpp{05}, the index  $$|B(x_1,  x_2, \dots , x_k;a):B(px_1,  x_2, \dots , x_k;a)|$$ is $p^{{\frac 1k} \sum_{a \in \mathcal A (n)} N(k,a)}$ =
$p^{{\frac 1k}\sum_{(a_j)=a   \in \mathcal A (n)} \prod_{j=1}^n a_j{{a_j+k-1}\choose {k-1}}}$.
\end{coro}

Now assume that  $H$ has a nondegenerate symmetric bilinear form $\<\cdot|\cdot\>.$ Then we can make $M(1)$ an irreducible module for the affine algebra
$$\hat H=H\otimes \CC[t,t^{-1}]\oplus \CC K$$
such that $H\otimes \CC[t]$  annihilates $\bf 1$ and the central element $K$ acts as $1.$   We abbrevuate  $h\otimes t^m$ by writing $h(m)$ for $h\in H$ and $m\in \ZZ.$

\begin{nota}\label{ndsbf}  There is a unique
nondegenerate symmetric bilinear form $\<\cdot|\cdot\>$ on $M(1)$ such that $\<\vac|\vac\>=1$ and $\<h(m)u|v\>=-\<u | h(-m)v\>$ for $u,v\in M(1)$ and $h\in H$ (see \cite{L}, \cite{DG}).
\end{nota}

Furthermore,  $B_L(n)=B(x_1,...,x_k; n)$ is a $\ZZ$-form of $M(1)_n.$  Note that if $L=\ZZ x_1+\cdots +\ZZ x_k$ is rational lattice of $H$ in the sense that for any $\alpha,\beta\in L$
$\<\alpha|\beta\>\in\QQ$, then  $B(x_1,...,x_k; n)$ is also a rational lattice, due to the form.
In the notation of Corollary \refpp{06}, we have

\begin{coro}\labtt{07}   Assume existence of  the form as in Notation \ref{ndsbf}.
For $k\ge 1$ and $n\ge 0$, define 
$$S(k,n):= {\frac 1k} \sum_{a \in \mathcal A (n)} N(k,a) ={\frac 1k}\sum_{a=(a_j)  \in \mathcal A (n)} \prod_{j=1}^n a_j{{a_j+k-1}\choose {k-1}}.$$
Then $S(k,0)=0$ and 
$$\det(B(px_1,  x_2, \dots , x_k; n)) = \det(B(x_1,  x_2, \dots , x_k; n))p^{2S(k,n)}$$
for all $k\ge 1$ and $n\ge 0$.    
\end{coro}

\begin{rem}\labtt {08} This presentation helps us understand the ``homogeneous  part''  of the standard integral form in the symmetric algebra spanned over $\ZZ$ by all monomials made from a basis. The integral form involves expressions like Schur functions which have fractional coefficients so are not in the homogeneous part.   We shall study the quotient of that integral form by its homogeneous part.
\end{rem}

\section{Integral forms of $M(1)$}

Let $L$ be a positive definite integral lattice with basis $x_1, \dots , x_k$ and we denote the form on $L$ by $\ip \cdot\cdot.$   We recall the standard integral form for the  lattice vertex operator algebra based on $L$.

 Note from \cite{flm} that $M(1): =\CC[x_i(-n)|i=1,...,k;  n>0]$ is the Heisenberg vertex operator algebra and $V_L=M(1)\otimes \CC^{\epsilon}[L]$ is the corresponding lattice vertex operator algebra where $\epsilon$ is a bimultiplicative map from $L\times L\to \<\pm 1\>$
such that $\epsilon(\alpha,\beta)\epsilon(\beta,\alpha)=(-1)^{\ip \alpha\beta}$ and $\epsilon(\alpha,\alpha)=(-1)^{\ip \alpha\alpha/2}$
and where $\CC^{\epsilon}[L]=\oplus_{\alpha\in L}\CC e^{\alpha}$ is the twisted group algebra. There is a unique nondegenerate symmetric invariant bilinear form $\ip\cdot\cdot$ on $V_L$ such that
$$\ip {e^{\alpha}}{e^{\beta}}=\delta_{\alpha+\beta,0}$$
and
$$\ip{\alpha(n)u}{v}=-\ip u{\alpha(-n)v}$$
 for all  $u,v\in V_L$ $\alpha\in L$ and $n\in {\ZZ}$ (see \cite{B1986}, \cite{DG}).

Recall the subring $A_L$ from Section 3. Then $A_L$ is a $\ZZ$-form of $M(1)$ in the sense that $A_L$ is a vertex algebra over $\Z,$ $\ip{u}{v}\in\Z$ for $u,v\in A_L$ and
$M(1)=\CC\otimes_{\Z}A_L$ \cite{DG}.

We also set $(V_L)_\ZZ=\oplus_{\alpha\in L}A_L\otimes e^{\alpha}.$  Then $(V_L)_\ZZ$ is a vertex operator algebra over $\ZZ$ generated by $e^{\pm x_i}$ for $i=1,...,k$ and is a free $\ZZ$-module
such that $V_L=\CC\otimes_{\ZZ}(V_L)_{\ZZ}.$

\begin{de}\label{standardintform} $(V_L)_\ZZ=\oplus_{\alpha\in L}A_L\otimes e^{\alpha}.$  is the {\it standard integral form}  in the lattice vertex operator algebra $V_L$.
\end{de}

Let  $(V_L)_{\ZZ,n}$  consists of vectors of weight $n$ in $(V_L)_{\ZZ}.$ To get $det((V_L)_{\ZZ, n})$, we first understand $\det(A_L(n))$ in terms of $\det(L).$   Recall
$B_{L}=\Z[x_i(-n)\mid i=1,...,k; \, n>0]$ and $B_L(n)=B_L\cap M(1)_n$ for $n\geq 0.$   By Lemma  \ref{l1}, $[A_L(n):B_L(n)]$ only depends on the rank of $L$ and integer $n.$

Using Lemma \ref{l1}  we can give  an explicit expression of  $[A_L(n):B_L(n)]$.   We define numbers $b_0:=1$ and for $n>0$,   
$$b_n : =\prod_{a=(a_1,a_2,\cdots)\in \mathcal A (n)}\prod_{i\geq 1}  i^{a_i}\cdot  a_i! .$$

\begin{lem}\label{l2} The index $[A_L(n):B_L(n)]$ is the square root of
$$\prod_{n_1,...,n_k\geq 0, \, \sum n_i=n }b_{n_1}\cdots b_{n_k}$$
 for $n\geq 0.$
\end{lem}

\pf By Lemma \ref{l1},   $[A_L(n):B_L(n)]$ is independent of lattice $L.$ So we can choose $L=\ZZ x_1+ \cdots +\ZZ x_k$ such that
 $\{x_1,...,x_k\}$ is an orthonormal basis of $H$ for convenience of computation.  Then
$A_L(n)$ is a unimodular lattice by Proposition 3.6 of \cite{DG}. It is easy to show that
\begin{eqnarray*}
& &\ \  \ \ip {x_i(-p)^s}{x_i(-p)^s}=(-1)^s\ip{{\bf 1}}{x_i(p)^sx_i(-p){\bf 1}}\\
& &=(-1)^ss!p^s\ip{{\bf 1}}{{\bf 1}}\\
& &=(-1)^ss!p^s
\end{eqnarray*}
for any $i,s$.    This shows that $|\det (B_L(n))|$ equals $\prod_{n_1,...,n_k\geq 0, \,  \sum n_i=n}b_{n_1}\cdots b_{n_k}.$ Since $|\det (B_L(n))|=[A_L(n):B_{L,n}]^2,$ the result follows immediately. \eop

\begin{lem}\label{l3} Let $A_1,A_2, C_1, C_2$ be lattices with the same rank such that $C_i \subset A_i$ for $i=1,2,$  $C_2\subset C_1.$  Then $|\det (A_2)|= \frac{[C_1:C_2]^2[A_1:C_1]^2|\det (A_1)}{[A_2:C_2]^2}.$
In particular, if $\det (A_1)=1$ and $[A_1:C_1]=[A_2:C_2]$ then $|\det (A_2)|=[C_1:C_2]^2.$
\end{lem}
\pf The result follows from the following relations
$$ |\det (C_i)|=|\det (A_i)|[A_i:C_i]^2,\det (C_2)=\det (C_1)[C_1:C_2]^2$$
for $i=1,2.$
\eop

\begin{thm}\label{t1} Let $L$ be an positive definite integral lattice with a base $\{x_1,...,x_k\}$ as before. Then  for each $n\geq 0,$ $|\det(A_L(n))|$ is an integer power of $\det(L).$   In fact,  $|\det (A_L(n))|=\det(L)^{2S(k,n)}$, where $S(k,n)$ is given by Lemma \refpp{07}.
\end{thm}
\pf We prove the theorem in several steps.   If $K$ contains the sublattice $J$ with finite index, then one may deduce the results for $K$ from those for $J$, and conversely, from the results of Section 3.

The result $|\det (A_L(n))|=1$ when $L$ is unimodular was proved in
\cite{DG}.   Let $p$ be a positive integer.

Case (a): Let $0\ne p \in \ZZ$ and  $L=\Z pe_1\oplus \Z e_2\oplus\cdots \oplus\Z e_k$  be a sublattice of $\Z^k=\Z e_1\oplus \cdots \oplus\Z e_k$ where $\{e_1,...,e_k\}$ is the standard orthonormal basis of $\RR^k.$
Using Lemmas \ref{l1}, \ref{l2} with $A_1=A_{\Z^k}(n),$ $A_2=A_L(n),$ $C_1=B_{\Z^k}(n)$ and $C_2=B_L(n)$ gives $|\det (A_L(n))|=[B_{\Z^k}(n):B_L(n)]^2.$ Note that $\det(L)=p^2.$
By Corollary \ref{07}, $|\det (A_L(n))|=p^{2S(k,n)}=\det (L)^{2S(k,n)}.$

Case (b):  Let  $L=\Z p_1e_1\oplus \Z p_2e_2\oplus\cdots \oplus\Z p_ke_k$  for any positive integers $p_1,...,p_k.$ Then $\det ( L)=p_1^2\cdots p_k^2$ and $|\det (A_L(n))|=\det (L)^{2S(k,n)}$ by Case (a).

Case (c): Let  $T$ be a positive integer such that $L$ is a rank $k$ sublattice of $K=\frac{1}{T}\Z^k$, i.e., $L \subset \QQ^k$.   Then $\det (L)=[K:L]^2T^{-2k}.$ There exist
a base $\{u_1,...,u_k\}$ of $K$ and positive integers $p_1,...,p_k$ such that $\{p_1u_1,...,p_ku_k\}$ is a base of $L.$ This implies that $[K:L]=p_1\cdots p_k.$ From Lemma \ref{l3} and
discussion in Case (b),
we see that $|\det (A_L(n))|=[B_K(n):B_L(n)]^2|\det (A_K(n))|=(p_1\cdots p_k)^{2S(k,n)}|\det (A_K(n))|.$ On the other hand,
$$1=|\det(A_{\Z^k}(n))|=[B_K(n): B_{\Z^k}(n)]^2|\det (A_K(n))|=T^{{2k S(k,n)}}|\det (A_K(n))|.$$ Thus
$$|\det (A_L(n))|=(p_1\cdots p_k)^{2S(k,n)} T^{-{2kS(k,n)}}=\det (L)^{2S(k,n)}.$$

Case (d): Let $L$ be an arbitrary integral lattice in Euclidean space $\RR^k.$
The problem with applying (c) is that $L$ is not necessarily a sublattice of $\QQ^k.$  However,  we can use a sequence of rational lattices $\{L_i|i\in\NN\}$ such that
``$\lim_{i\to \infty}L_i=L$''.    We fix a base $\{v_1,...,v_k\}$ of $L$ and take linearly independent vectors $v^i_1,...,v^i_k\in \QQ^k$ such
that $|v^i_j-v_j|<\frac{1}{i}$ for all $i,j.$   It is clear that $\lim_{i\to \infty}\det (L_i)=\det (L)$ and $\lim_{i\to\infty}|\det (A_{L_i}(n))|=|\det  (A_L(n))|$ for any $n\geq 0.$
It follows from Case (c) that $|\det (A_L(n))|=\det (L)^{2S(k,n)},$ as desired. \eop

\section{Integral forms of $V_L$}

We now assume that $L$ is a positive definite even lattice. Recall that  $(V_L)_\ZZ=\oplus_{\alpha\in L}A_L\otimes e^{\alpha}.$ Also recall $(V_L)_{\ZZ,n}$ from Section 4.  We determine $\det ((V_L)_{\ZZ,n})$ in this section.

For $m\geq 0$ we set $L_{2m}=\{\alpha\in L|\ip \alpha\alpha=2m\}.$   Define $Y_0:=L_o=\{0\}$.    
For $m\ge 1$, let $Y_{2m}$ be a subset of $L_{2m}$ such that $2|Y_{2m}|=|L_{2m}|$ and $L_{2m}=Y_{2m}\cup (-Y_{2m}).$  For $\alpha\in L$ we set
$W^{\alpha}=M(1)_{L}\otimes e^{\alpha}+M(1)_L\otimes e^{-\alpha}\subset (V_L)_{\ZZ}.$ Let $W^{\alpha}_n=W^{\alpha}\cap (V)_{L,n}.$ Then  $W^{\alpha}_n\ne 0$ if and only if $\alpha\in L_{2m}$ and $m\leq n.$
In this case, $W^{\alpha}_n=A_L(n-m)\otimes e^{\alpha}+A_L(n-m)\otimes e^{-\alpha}.$   Observe that

$$(V_L)_{\ZZ,n}=\oplus_{m=0}^n\oplus_{\alpha\in Y_{2m}}W^{\alpha}_n$$
and $\ip {W^{\alpha}}{W^{\beta}}=0$ if $\alpha\ne \beta.$ So
$$\det( (V_L)_{\ZZ,n})=\prod_{m=0}^n\prod_{\alpha\in Y_{2m}}\det (W^\alpha_m).$$
From the definition of the bilinear form,   we know that
for $\alpha\in Y_{2m}$ with $m\ne 0$
$$|\det (W^{\alpha}_n)|=\det (A_L(n-m))^2.$$
Also, $det(W_0^0)=det(\ZZ \vac)=1$.   

Here is our main theorem, an immediate consequence of Theorem \ref{t1}.

\begin{thm}\label{maintheorem}  For all $n\ge 0$, we have
$$|\det ((V_L)_{\ZZ,n})|=\prod_{m=0}^n\det( L)^{|L_{2m}|{2S(k,n-m)}}.$$
\end{thm}

\section{Acknowledgements} 

C. Dong was supported by China NSF grant 11871351.
R. Griess was supported by funds from his Collegiate Professorship and Distinguished University Professorship at the University of Michigan.


\begin{thebibliography}{GRC99}

\bibitem{B1986}  R.  Borcherds,  Vertex algebras, Kac-Moody algebras, and the Monster, {\em Proc. Nat. Acad. Sci. U.S.A.} {\bf 83} (1986),
    3068-3071.
\bibitem{DG} C.  Dong and R. L. Griess Jr, Integral  forms in vertex operator algebras which are invariant under finite groups, {\em J.
Algebra} {\bf 365} (2012), 184-198.
\bibitem{DG2} C.  Dong and R. L. Griess Jr,  Lattice-integrality of certain group-invariant integral forms in vertex operator algebras, {\em J. Algebra} {\bf 474} (2017), 505-516.

\bibitem{flm} I. B. Frenkel, J. Lepowsky and A. Meurman,
Vertex Operator Algebras and the Monster, {\em Pure and Applied
Math.,} Vol. {\bf 134}, Academic Press, 1988.

\bibitem{L} H. Li, Symmetric invariant bilinear forms on vertex operator
algebras, {\em Pure and Appl. Math. } {\bf 96} (1994), 279-297.
\end{thebibliography}
\end{document}